\documentclass{article}
\usepackage{amsmath,amsthm,amsopn,amstext,amscd,amsfonts,amssymb}
\usepackage{natbib}
\usepackage{verbatim}

\def\tr{\mathop{\rm tr}\nolimits}
\def\R{\mathop{\rm Re}\nolimits}

\def\etr{\mathop{\rm etr}\nolimits}

\renewenvironment{abstract}
                 {\vspace{6pt}
                  \begin{center}
                  \begin{minipage}{5in}
                  \centerline{\textbf{Abstract}}
                  \noindent\ignorespaces
                 }
                 {\end{minipage}\end{center}}
\newtheorem{thm}{\textbf{Theorem}}[section]
\newtheorem{cor}{\textbf{Corollary}}[section]

\theoremstyle{definition}

\title{\huge \textbf{Noncentral bimatrix variate generalised beta distributions}}
\author{
  \textbf{Jos\'e A. D\'{\i}az-Garc\'{\i}a} \thanks{Corresponding author\newline
    }\\
  Department of Statistics and Computation \\
  25350 Buenavista, Saltillo, Coahuila, Mexico \\
  E-mail: jadiaz@uaaan.mx \\[2ex]
  \textbf{Ram\'on Guti\'errez J\'aimez} \\
  Department of Statistics and O.R. \\
  University of Granada \\
  Granada 18071, Spain \\
  E-mail: rgjaimez@ugr.es\\
}
\date{}
\begin{document}
\maketitle

\begin{abstract}
In this paper, we determine the density functions of nonsymmetrised doubly noncentral
matrix variate beta type I and II distributions. The nonsymetrised density functions
of doubly noncentral and noncentral bimatrix variate generalised beta type I and II
distributions are also obtained.
\end{abstract}

\section{Introduction}\label{intro}

When we consider generalising the distribution of a random variable to the
multivariate case, two options are normally addressed, those of extending it to
either the vectorial or the matrix cases, e.g. normal, t or bessel distributions,
among many others. However, some of these generalisations have traditionally been
made directly to the matrix case, where such a matrix is symmetric - this is the case
of the chi-square and beta distributions, for which the corresponding multivariate
distributions are the Wishart and matrix variate beta distributions, respectively.
Nevertheless, these latter generalisations are inappropriate in some cases, because
sometimes we might be interested in a vectorial version and not in a matrix version.
For example, we are interested in a random vector in which each marginal is a random
variable beta (type I or II). \citet{ln:82} proposed a multivariate (vector) beta
distribution. Some applications to utility modelling and Bayesian analysis are also
presented in \citet{ln:82} and \citet{cn:84}, respectively. In particular
\citet{ol:03} proposed the following bivariate version. Observe that the following
definition eliminates the hypothesis that the variables have a chi-squared
distribution, assuming, instead, a gamma distribution.

Let $A$, $B$ and $C$ be distributed as independent gamma random variables with
parameters $\alpha = a, b, c$, respectively and $\delta = 1$ in the three cases (see
eq. (\ref{gammac}) in Section \ref{sec2}), and define
\begin{equation}\label{eq1}
    U_{1} = \frac{A}{A+C}, \qquad U_{2} = \frac{B}{B+C}.
\end{equation}
Clearly, $U_{1}$ and $U_{2}$ each have a beta type I distribution, $U_{1} \sim
\mathcal{B}I_{1}(a,c)$ and $U_{2} \sim \mathcal{B}I_{1}(b,c)$, over $0 \leq
u_{1},u_{2} \leq 1$. However, they are correlated, and then $(U_{1},U_{2})'$ has a
bivariate generalised beta type I distribution over $0 \leq u_{1},u_{2} \leq 1$.

A similar result is obtained in the case of beta type II. Now, let us define
$$
  F_{1} = \frac{A}{C}, \qquad F_{2} = \frac{B}{C}.
$$
Once again it is easy to see that $F_{1}$ and $F_{2}$ each have a beta type II
distribution, $F_{1} \sim \mathcal{B}II_{1}(a,c)$ and $F_{2} \sim
\mathcal{B}II_{1}(b,c)$, over $f_{1},f_{2} \geq 0$. As in the beta type I case, they
are correlated therefore and $(F_{1},F_{2})'$ has a bivariate generalised beta type
II distribution over $f_{1},f_{2} \geq 0$.

These ideas can be extended to the matrix variate case. Thus, let us assume a
partitioned matrix $\mathbb{U} = (\mathbf{U}_{1} \vdots \mathbf{U}_{2})' \in \Re^{2m
\times m}$, then under the matrix variate versions of the transformations
(\ref{eq1}), we are interested in finding the joint density of $\mathbf{U}_{1}$ and
$\mathbf{U}_{2}$, where it is easy to see that the marginal densities of
$\mathbf{U}_{1}$ and $\mathbf{U}_{2}$ are matrix variate beta type I distributions.
In the central case, the matrix variate joint densities of $\mathbf{U}_{1}$ and
$\mathbf{U}_{2}$ and of $\mathbf{F}_{1}$ and $\mathbf{F}_{2}$ and some properties are
studied in \citet{dggj:08b}. These distributions are termed central bimatrix variate
generalised beta type I and II distributions, respectively. They play a potentially
important role in the context of shape theory, specifically in affine or
configuration densities, such as the \citet{gm:93} conjecture. Suppose that we have
two samples of images of size $n$, each one of which is obtained at two times. Also,
assume that we are interested in evaluating whether a learning process is present or
whether the process has a memory. In this context, if we obtain as the configuration
density a central, noncentral or double noncentral bimatrix variate generalised beta
type I and II distribution,  it might be possible to study these problems (learning
or memory problems) and to compare the parameters of $\mathbf{F}_{1}$
($\mathbf{U}_{1}$) and $\mathbf{F}_{2}$ ($\mathbf{U}_{2}$) considering the latter as
bimatrix variate.

In this paper, we study bimatrix variate generalised beta type I and II distributions
under different cases of noncentrality. Some definitions regarding the symmetrised
function are given in Subsection \ref{ssec21} and Subsection \ref{ssec22} presents
known and new results about central, noncentral and doubly noncentral matrix variate
beta type I and II distributions; also we include the definition of the central
bimatrix variate generalised beta type I and II distributions.  Nonsymetrised doubly
noncentral density functions of the bimatrix variate generalised beta type I and II
distributions are studied, and diverse noncentral cases of the bimatrix variate
generalised beta type I and II distributions are obtained as particular cases of
nonsymetrised doubly noncentral density functions, see Sections \ref{sec3} and
\ref{sec4}, respectively.

\section{Preliminary results}\label{sec2}

\subsection{Symmetrised density function}\label{ssec21}

In multivariate analysis there exist a large class of important hypothesis testing
problems all of which may be tested by a set of criteria that depend functionally on
the eigenvalues of a matrix variate. With the propose to investigate the non-null
distributions of these criteria, \citet{g:73} introduce the notion of a
\textbf{symmetrised distribution} of a matrix variate, a notion which facilitates
many proofs in such derivations.

Given a density function $f_{_{\mathbf{X}}}(\mathbf{X})$, $\mathbf{X} \in \Re^{m
\times m}$, $\mathbf{X} > \mathbf{0}$, \citet{g:73} proposes the following definition
$$
  f_{s}(\mathbf{X}) = \int_{\mathcal{O}(m)} f_{_{\mathbf{X}}}(\mathbf{HXH}')
  (d\mathbf{H}), \quad \mathbf{H} \in \mathcal{O}(m)
$$
where $\mathcal{O}(m) = \{\mathbf{H}\in \Re^{m \times m}| \mathbf{H}'\mathbf{H} =
\mathbf{HH}' = \mathbf{I}_{m}\}$ and $(d\mathbf{H})$ denotes the normalised Haar
measure on $\mathcal{O}(m)$, see \citet[pp. 60 and 260]{mh:82}. This function
$f_{s}(\mathbf{X})$ is termed symmetrised density function of $\mathbf{X}$.

Our proposal is to apply this idea from \citet{g:73} in an inverse way, i.e.
well-known the explicit expression of the symmetrised density function of X
\begin{equation}\label{ne}
    f_{s}(\mathbf{X}) = \int_{\mathcal{O}(m)}f(\mathbf{HXH}')(d\mathbf{H}).
\end{equation}
We wish to identify the density function $f(\mathbf{X})$. The density function
obtained by applying the idea underlying (\ref{ne}) is termed the nonsymmetrised
density function. Finally, note that the joint density function of the eigenvalues of
$\mathbf{X}$ can be found from $f_{s}(\mathbf{X})$ or $f(\mathbf{X})$, indifferently.

\subsection{Matrix variate beta distributions}\label{ssec22}

In general, matrix variate beta type I and II distributions are defined in terms of
two matrices, say, $\mathbf{A}$ and $\mathbf{B}$, which are independent and have
Wishart distributions, see \citet{or:64}, \citet{k:70}, \citet{mh:82}, \citet{f:85},
\citet{c:96}, \citet{gn:00}, \citet{dggj:07, dggj:08a, dggj:08b}, among many others.
The present paper generalises these results, assuming that $\mathbf{A}$ and
$\mathbf{B}$ have matrix variate gamma distributions.

The $m \times m$ matrix $\mathbf{A}$ is said to have a noncentral matrix variate
gamma distribution with parameters $a \in \Re$ in which $\mathbf{\Theta}$ is an $m
\times m$ positive definite matrix and $\mathbf{\Omega}$ is an $m \times m$ matrix,
this fact being denoted as $\mathbf{A} \sim \mathcal{G}_{m}(a, \mathbf{\Theta},
\mathbf{\Omega})$, if its density function is (see \citet[pp. 57 and 61]{mh:82} and
\citet{gn:00})
\begin{equation}\label{gamma}
    \mathcal{G}_{m}(\mathbf{A};a,\mathbf{\Theta},\mathbf{\Omega}) =
    \mathcal{G}_{m}(\mathbf{A};a,\mathbf{\Theta})
    {}_{0}F_{1}(a,\mathbf{\Omega \Theta}^{-1}\mathbf{A}), \quad \mathbf{A} > \mathbf{0},
\end{equation}
where ${}_{0}F_{1}(\cdot)$ is a hypergeometric function with a matrix argument (see
\citet[p. 258]{mh:82}) and $\mathcal{G}_{m}(\mathbf{A};a,\mathbf{\Theta}) \equiv
\mathcal{G}_{m}(\mathbf{A};a,\mathbf{\Theta},\mathbf{0})$ denotes the density
function of a central matrix variate gamma distribution given by
\begin{equation}\label{gammac}
    \mathcal{G}_{m}(\mathbf{A};a,\mathbf{\Theta}) =
    \frac{|\mathbf{A}|^{a-(m+1)/2}}{\Gamma_{m}[a]|\mathbf{\Theta}|^{a}}\etr(-\mathbf{\Theta}^{-1}\mathbf{A}),
    \quad \mathbf{A} > \mathbf{0},
\end{equation}
and denoted as $\mathbf{A} \sim \mathcal{G}_{m}(a, \mathbf{\Theta}) \equiv
\mathcal{G}_{m}(a, \mathbf{\Theta}, \mathbf{0})$. Where $\etr(\cdot) \equiv
\exp(\tr(\cdot))$ and $\Gamma_{m}[a]$ denotes the multivariate gamma function and is
defined as
$$
  \Gamma_{m}[a] = \int_{\mathbf{V}>0} \etr(-\mathbf{V}) |\mathbf{V}|^{a-(m+1)/2} (d\mathbf{V}),
$$
$\R(a) > (m-1)/2$.

In addition to the classification of beta distributions as beta type I and type II
(see \citet{gn:00} and \citet{sk:79}), two definitions have been proposed for each of
these, see \citet{or:64}, \citet{s:68}, \citet{dggj:01} and \citet{j:64}. Let us
focus initially on the beta type I distribution; if $\mathbf{A}$ and $\mathbf{B}$
have a matrix variate gamma distribution, i.e. $\mathbf{A} \sim
\mathcal{G}_{m}(a,\mathbf{I}_{m})$ and $\mathbf{B} \sim
\mathcal{G}_{m}(b,\mathbf{I}_{m})$ independently, then the beta matrix $\mathbf{U}$
can be defined as
\begin{equation}\label{defbI}
  \mathbf{U} =
  \left\{%
     \begin{array}{ll}
      (\mathbf{A} + \mathbf{B})^{-1/2}\mathbf{A} ((\mathbf{A} + \mathbf{B})^{-1/2})', & \mbox{Definition $1$ or},\\
       \mathbf{A} ^{1/2}(\mathbf{A} + \mathbf{B})^{-1} (\mathbf{A} ^{1/2})', & \mbox{Definition $2$},\\
    \end{array}%
  \right.
\end{equation}
where $\mathbf{C}^{1/2}(\mathbf{C}^{1/2})' = \mathbf{C}$ is a reasonable nonsingular
factorization of $\mathbf{C}$, see \citet{gn:00}, \citet{sk:79} and \citet{mh:82}. It
is readily apparent that under definition $1$ and $2$ the density function is
\begin{equation}\label{beta}\hspace{-0.7cm}
    \mathcal{B}I_{m}(\mathbf{U};a,b) =
    \frac{1}{\beta_{m}[a,b]} |\mathbf{U}|^{a -(m + 1)/2} |\mathbf{I}_{m} - \mathbf{U}|^{b -(m +1)/2}
    (d\mathbf{U}), \quad \mathbf{0} < \mathbf{U} < \mathbf{I}_{m},
\end{equation}
writing this fact as $\mathbf{U} \sim \mathcal{B}I_{m}(a,b)$, with $\R(a) > (m-1)/2$
and $\R(b) > (m-1)/2$; where $\beta_{m}[a,b]$ denotes the multivariate beta function
defined by
\begin{eqnarray*}
  \beta_{m}[b,a] &=& \int_{\mathbf{0}<\mathbf{S}<\mathbf{I}_{m}} |\mathbf{S}|^{a-(m+1)/2}
  |\mathbf{I}_{m} - \mathbf{S}|^{b-(m+1)/2} (d\mathbf{S}) \\
    &=& \int_{\mathbf{R}>\mathbf{0}} |\mathbf{R}|^{a-(m+1)/2} |\mathbf{I}_{m} + \mathbf{R}|^{-(a+b)} (d\mathbf{R}) \\
    &=& \frac{\Gamma_{m}[a] \Gamma_{m}[b]}{\Gamma_{m}[a+b]}.
\end{eqnarray*}

A similar situation arises with the beta type II distribution, and thus we have the
following two definitions:
\begin{equation}\label{defbII}
  \mathbf{F} =
  \left\{%
     \begin{array}{ll}
      \mathbf{B}^{-1/2}\mathbf{A} (\mathbf{B}^{-1/2})', & \mbox{Definition 1},\\
      \mathbf{A}^{1/2}\mathbf{B}^{-1} (\mathbf{A} ^{1/2})', & \mbox{Definition 2},\\
    \end{array}%
  \right.
\end{equation}
with the distribution being denoted as $\mathbf{F} \sim \mathcal{B}II_{m}(a,b)$. In
this case, under definitions 1 and 2, the density function of $\mathbf{F}$ is
\begin{equation}\label{efe}
    \mathcal{B}II_{m}(\mathbf{F};a,b) =
    \frac{1}{\beta_{m}[a,b]}  |\mathbf{F}|^{a-(m+1)/2}|\mathbf{I}_{m} +
    \mathbf{F}|^{-(a+b)}, \quad \mathbf{F} > 0.
\end{equation}

\citet{dggj:07,dggj:08a} showed that in doubly noncentral and noncentral matrix
variate beta type I and II distributions, the corresponding  density functions are
invariant under definitions 1 and 2. Therefore, henceforth we shall make no
distinction between definitions 1 and 2.

When these ideas are extended to the doubly noncentral case, i.e. when $\mathbf{A}
\sim \mathcal{G}_{m}(a,\mathbf{I}_{m}, \mathbf{\Omega}_{1})$ and $\mathbf{B} \sim
\mathcal{G}_{m}(b,\mathbf{I}_{m}, \mathbf{\Omega}_{2})$, strictly speaking, we have
not found the densities of the matrix variate beta type I and II distributions.
Rather, for the case of the beta type II distribution, \citep{ch:80} found the
distribution of $\tilde{\mathbf{V}} =
\tilde{\mathbf{B}}^{-1/2}\tilde{\mathbf{A}}(\tilde{\mathbf{B}}^{-1/2})'$ where
$\tilde{\mathbf{A}} = \mathbf{H}'\mathbf{AH}$ and $\tilde{\mathbf{B}} =
\mathbf{H}'\mathbf{BH}$, $\mathbf{H} \in \mathcal{O}(m)$, with $\mathcal{O}(m) =
\{\mathbf{H} \in \Re^{m \times m}| \mathbf{HH}' = \mathbf{H}'\mathbf{H} =
\mathbf{I}_{m}\}$. It is straightforward to show that the procedure proposed by
\citet{ch:80} and \citet{chd:86} is equivalent to finding the symmetrised density
defined by \citet{g:73}, see also \citet{r:75}. From \citet{dggj:08a} and using the
notation for the operator sum as in \citet{da:80} we have the following:
\begin{enumerate}
  \item The symmetrised density function of doubly noncentral matrix variate beta type
  I is
  \begin{equation}\label{sdnbt1}
    \mathcal{B}I_{m}(\mathbf{U};a,b) \etr\left(-(\mathbf{\Omega}_{1}+ \mathbf{\Omega}_{2})\right)
    \hspace{6cm}
  \end{equation}
  \par\indent\hfill\mbox{$\displaystyle{\times \ \
    \sum_{\kappa,\lambda; \ \phi}^{\infty}\frac{(a+b)_{\phi}}{(a)_{\kappa} (b)_{\lambda} k! \ l!}
    \frac{C_{\phi}^{\kappa,\lambda}(\mathbf{\Omega}_{1}, \mathbf{\Omega}_{2}) C_{\phi}^{\kappa,\lambda}
    (\mathbf{U}, (\mathbf{I}_{m}-\mathbf{U}))}{C_{\phi}(\mathbf{I}_{m})}}, \ \ \mathbf{0} < \mathbf{U} <\mathbf{I}_{m}$.}
  \par\noindent
  \item and the symmetrised density function of doubly noncentral matrix variate beta type
  II is
  \begin{equation}\label{sdnbt2}
    \mathcal{B}II_{m}(\mathbf{F};a,b) \etr\left(-(\mathbf{\Omega}_{1}+ \mathbf{\Omega}_{2})\right)
    \hspace{5cm}
  \end{equation}
  \par\indent\hfill\mbox{$\displaystyle{\times \ \
    \sum_{\kappa,\lambda;\ \phi}^{\infty}\frac{(a+b)_{\phi}}{(a)_{\kappa} \ (b)_{\lambda} k! \ l!}
    \frac{C_{\phi}^{\kappa,\lambda}(\mathbf{\Omega}_{1}, \mathbf{\Omega}_{2}) C_{\phi}^{\kappa,\lambda}
    ((\mathbf{I}_{m} + \mathbf{F})^{-1} \mathbf{F}, (\mathbf{I}_{m} + \mathbf{F})^{-1})}
    {C_{\phi}(\mathbf{I}_{m})}}.$}
  \par\noindent
\end{enumerate}
where $\mathbf{F} > \mathbf{0}$, $\R(a) > (m-1)/2$, $\R(b) > (m-1)/2$, $(a)_{\tau}$
is the generalised hypergeometric coefficient or product of Poch\-hammer symbols and
$C_{\phi}^{\kappa,\lambda}(\cdot,\cdot)$ denotes the invariant polynomials with the
matrix arguments defined in \citet{da:80}, see also \citet{ch:80} and \cite{chd:86}.

As particular cases of doubly noncentral distributions it is possible to obtain two
different definitions of noncentral distributions, given another classification,  in
which the beta matrix is defined as follows, see \citet{g:73} and \citet{gn:00}:
\begin{eqnarray*}\hspace{-.5cm}
  \mathbf{W} &=& \mathbf{A} ^{1/2}(\mathbf{A} + \mathbf{B})^{-1} (\mathbf{A} ^{1/2})',
    \ \mbox{ denoting as } \ \mathcal{B}I(A)_{m}(a,b,\mathbf{\Omega}_{2}),
    \ (\mathbf{\Omega}_{1} = \mathbf{0})\\
  \mathbf{U} &=& \mathbf{A} ^{1/2}(\mathbf{A} + \mathbf{B})^{-1} (\mathbf{A} ^{1/2})',
    \ \mbox{ denoting as } \ \mathcal{B}I(B)_{m}(a,b,\mathbf{\Omega}_{1}),
    \ (\mathbf{\Omega}_{2} = \mathbf{0}).
\end{eqnarray*}
Similarly, in the case of beta type II we have
\begin{eqnarray*}
  \mathbf{V} &=& \mathbf{B}^{-1/2}\mathbf{A} (\mathbf{B}^{-1/2})', \quad \mbox{denoting as}
    \quad \mathcal{B}II(A)_{m}(a,b,\mathbf{\Omega}_{2}),
    \ (\mathbf{\Omega}_{1} = \mathbf{0})\\
  \mathbf{F} &=& \mathbf{B}^{-1/2}\mathbf{A} (\mathbf{B}^{-1/2})', \quad \mbox{denoting as}
  \quad \mathcal{B}II(B)_{m}(a,b,\mathbf{\Omega}_{1}),
    \ (\mathbf{\Omega}_{2} = \mathbf{0})
\end{eqnarray*}
Both distributions, types A and B, play a fundamental role in various areas of
statistics, for example in the $W$ and $U$ criteria proposed by \citet{dgcl:08}.

The symmetrised and nonsymmetrised density functions of $\mathbf{W}$, $\mathbf{U}$,
$\mathbf{V}$ and $\mathbf{F}$ can be obtained as particular cases of (\ref{sdnbt1})
and (\ref{sdnbt2}). All these densities are found in \citet{dggj:07}.

Now, using the approach described in \citet{dggj:07} we can find the (nonsymmetrised)
density functions of doubly noncentral matrix variate beta type I and II
distributions. Observe that for (\ref{sdnbt1}) and by \citet{dg:06},
$$
   \int_{\mathcal{O}(m)} C_{\phi}^{\kappa,\lambda}
    (\mathbf{\Omega}_{1}\mathbf{HUH}', \mathbf{\Omega}_{2}(\mathbf{I}_{m}-\mathbf{HUH}')) (d\mathbf{H})=
    \frac{C_{\phi}^{\kappa,\lambda}(\mathbf{\Omega}_{1}, \mathbf{\Omega}_{2}) C_{\phi}^{\kappa,\lambda}
    (\mathbf{U}, (\mathbf{I}_{m}-\mathbf{U}))}{\theta_{\phi}^{\kappa,\lambda} \
    C_{\phi}(\mathbf{I}_{m})}.
$$
where $\theta_{\phi}^{\kappa,\lambda}$ is defined in \citet{da:79} and \citet{ch:80}.

Proceeding in analogous form for (\ref{sdnbt2}), we have the following.
\begin{thm}\label{teo1}
For $\R(a)> (m-1)/2$ and $\R(b)> (m-1)/2$,
\begin{enumerate}
  \item the nonsymmetrised density function of the doubly noncentral matrix variate beta type
  I is
  \begin{equation}\label{dnbt1}
    \mathcal{B}I_{m}(\mathbf{U};a,b) \etr\left(-(\mathbf{\Omega}_{1}+ \mathbf{\Omega}_{2})\right)
    \hspace{6cm}
  \end{equation}
  \par\indent\hfill\mbox{$\displaystyle{\times \ \
    \sum_{\kappa,\lambda; \ \phi}^{\infty}\frac{(a+b)_{\phi}\ \theta_{\phi}^{\kappa,\lambda}}{(a)_{\kappa}
    (b)_{\lambda} k! \ l!}  \ C_{\phi}^{\kappa,\lambda}\left(\mathbf{\Omega}_{1}\mathbf{U},
    \mathbf{\Omega}_{2}(\mathbf{I}_{m}-\mathbf{U})\right),}
    \ \ \mathbf{0} < \mathbf{U} <\mathbf{I}_{m}$.}
  \par\noindent
  which is denoted as $\mathbf{U} \sim \mathcal{B}I_{m}(a, b, \mathbf{\Omega}_{1}, \mathbf{\Omega}_{2})$.
  \item and the nonsymmetrised density function of the doubly noncentral matrix variate beta type
  II is
  \begin{equation}\label{dnbt2}
    \mathcal{B}II_{m}(\mathbf{F};a,b) \etr\left(-(\mathbf{\Omega}_{1}+ \mathbf{\Omega}_{2})\right)
    \hspace{5cm}
  \end{equation}
  \par\indent\hfill\mbox{$\displaystyle{\times \ \
    \sum_{\kappa,\lambda; \ \phi}^{\infty}\frac{(a+b)_{\phi}\ \theta_{\phi}^{\kappa,\lambda}}{(a)_{\kappa}
    \ (b)_{\lambda} k! \ l!} \ C_{\phi}^{\kappa,\lambda} (\mathbf{\Omega}_{1}(\mathbf{I}_{m} +
    \mathbf{F})^{-1} \mathbf{F}, \mathbf{\Omega}_{2}(\mathbf{I}_{m} + \mathbf{F})^{-1})}, \ \
    \mathbf{F} > \mathbf{0},$}
  \par\noindent
  which is denoted as $\mathbf{U} \sim  \mathcal{B}II_{m}(a, b, \mathbf{\Omega}_{1}, \mathbf{\Omega}_{2}).$
\end{enumerate}
\end{thm}

Doubly noncentral, noncentral and central matrix variate beta type I and II
distributions play a very important role in diverse problems  for proving hypotheses
in the context of multivariate analysis, including canonical correlation analysis,
the general linear hypothesis in MANOVA and multiple matrix variate correlation
analysis, see \citet{mh:82}, \citet{r:73}, \citet{s:68} and \citet{k:61}. Similarly,
doubly noncentral and noncentral beta distributions are to be found in the context of
econometrics and shape theory, see \citet{chd:86} and \citet{gm:93}, respectively.

Now from \citet{dggj:08b}; let $\mathbf{A}$, $\mathbf{B}$ and $\mathbf{C}$ be
independent, where $\mathbf{A} \sim \mathcal{G}_{m}(a, \mathbf{I}_{m})$, $\mathbf{B}
\sim \mathcal{G}_{m}(b, \mathbf{I}_{m})$ and $\mathbf{C} \sim \mathcal{G}_{m}(c,
\mathbf{I}_{m})$ with $\R(a)
> (m-1)/2$, $\R(b) > (m-1)/2$ and $\R(c) > (m-1)/2$ and let us define
\begin{equation}\label{bgb1}\hspace{-.7cm}
    \mathbf{U}_{1} = (\mathbf{A}+\mathbf{C})^{-1/2}\mathbf{A}(\mathbf{A} + \mathbf{C})^{-1/2}
  \quad \mbox{and} \quad \mathbf{U}_{2} = (\mathbf{B}+\mathbf{C})^{-1/2}\mathbf{B}
  (\mathbf{B}+\mathbf{C})^{-1/2}
\end{equation}
Of course, $\mathbf{U}_{1} \sim \mathcal{B}I_{m}(a,c)$ and $\mathbf{U}_{2} \sim
\mathcal{B}I_{m}(b,c)$. However, they are correlated and therefore the distribution
of $\mathbb{U} = (\mathbf{U}_{1}\vdots \mathbf{U}_{2})' \in \Re^{2m \times m}$ is
termed a central bimatrix variate generalised beta type I distribution, denoted as
$\mathbb{U} \sim \mathcal{BGB}I_{2m \times m}(a,b,c)$. Moreover, its density function
is
\begin{equation}\label{density1}\hspace{-0.7cm}
  \frac{|\mathbf{U}_{1}|^{a-(m+1)/2} |\mathbf{U}_{2}|^{b-(m+1)/2} |\mathbf{I}_{m} - \mathbf{U}_{1}|
  ^{b + c -(m+1)/2} |\mathbf{I}_{m} - \mathbf{U}_{2}|^{a + c -(m+1)/2}}{\beta^{*}_{m}[a,b,c]
  |\mathbf{I}_{m} - \mathbf{U}_{1}\mathbf{U}_{2}|^{a + b + c }}
\end{equation}
and is denoted as $\mathcal{BGB}I_{2m \times m}(\mathbb{U};a,b,c)$, where $\mathbf{0}
< \mathbf{U}_{1}< \mathbf{I}_{m}$, $\mathbf{0} < \mathbf{U}_{2}< \mathbf{I}_{m}$ with
$\R(a)> (m-1)/2$, $\R(b) > (m-1)/2$ and $\R(c)> (m-1)/2$ and
$$
  \beta^{*}_{m}[a,b,c] = \frac{\Gamma_{m}[a] \Gamma_{m}[b] \Gamma_{m}[c]}{\Gamma_{m}[a+b+c]}
$$

Similarly, let
\begin{equation}\label{bgb2}
    \mathbf{F}_{1} = \mathbf{C}^{-1/2}\mathbf{A}\mathbf{C}^{-1/2}
  \quad \mbox{and} \quad \mathbf{F}_{2} = \mathbf{C}^{-1/2}\mathbf{B}
  \mathbf{C}^{-1/2}
\end{equation}
Clearly, $\mathbf{F}_{1} \sim \mathcal{B}II_{m}(a,c)$ and $\mathbf{F}_{2} \sim
\mathcal{B}II_{m}(b,c)$. But they are correlated and then the distribution of
$\mathbb{F} = (\mathbf{F}_{1}\vdots \mathbf{F}_{2})' \in \Re^{2m \times m}$ is termed
a central bimatrix variate generalised beta type II distribution, which is denoted as
$\mathbb{F} \sim \mathcal{BGB}II_{2m \times m}(a,b,c)$. Its density function is
\begin{equation}\label{density2}
    \mathcal{BGB}II_{2m \times m}(\mathbb{F};a,b,c) =
    \frac{|\mathbf{F}_{1}|^{a-(m+1)/2} |\mathbf{F}_{2}|^{b-(m+1)/2} }{\beta^{*}_{m}[a,b,c]
  |\mathbf{I}_{m} + \mathbf{F}_{1} + \mathbf{F}_{2}|^{a + b + c }}
\end{equation}
where $\mathbf{F}_{1} > \mathbf{0}$, $\mathbf{F}_{2} > \mathbf{0}$, with $\R(a)>
(m-1)/2$, $\R(b) > (m-1)/2$ and $\R(c) > (m-1)/2$.

Other properties of bimatrix variate generalised beta type I and II distributions are
studied in \cite{dggj:08b}.

The use of matrix and bimatrix variate beta-type distributions has not been developed
as expected and hoped, due particularly to the fact that such distributions depend on
hypergeometric functions with a matrix argument or on zonal polynomials, which until
very recently were quite complicated to evaluate. Recently, descriptions have been
made of algorithms that are very efficient at calculating both zonal polynomials and
hypergeometric functions with a matrix argument; these can be used more widely and
more efficiently in noncentral distributions in general, see \citet{k:04} and
\citet{ke:06}.

\section{Doubly noncentral bimatrix variate generalised beta type I distribution}\label{sec3}

In this section we derive the doubly noncentral bimatrix variate generalised beta
type I distribution.

\begin{thm}\label{dngb1}
Let $\mathbf{A}$, $\mathbf{B}$ and $\mathbf{C}$ be independent random matrices, such
that $\mathbf{A} \sim \mathcal{G}_{m}(a, \mathbf{I}_{m}, \mathbf{\Omega}_{1})$,
$\mathbf{B} \sim \mathcal{G}_{m}(b, \mathbf{I}_{m}, \mathbf{\Omega}_{2})$ and
$\mathbf{C} \sim \mathcal{G}_{m}(c, \mathbf{I}_{m}, \mathbf{\Omega}_{3})$ with
$\R(a)> (m-1)/2$, $\R(b) > (m-1)/2$ and $\R(c) > (m-1)/2$ and let us define
\begin{equation}\label{bgb12}\hspace{-.7cm}
    \mathbf{U}_{1} = (\mathbf{A}+\mathbf{C})^{-1/2}\mathbf{A}(\mathbf{A} + \mathbf{C})^{-1/2}
  \quad \mbox{and} \quad \mathbf{U}_{2} = (\mathbf{B}+\mathbf{C})^{-1/2}\mathbf{B}
  (\mathbf{B}+\mathbf{C})^{-1/2}
\end{equation}
Then the symmetrised density function of  $\mathbb{U} = (\mathbf{U}_{1}\vdots
\mathbf{U}_{2})' \in \Re^{2m \times m}$ is
\begin{equation}\label{sdnbgbt1}\hspace{-0.7cm}
    \mathcal{BGB}I_{2m \times m}(\mathbb{U};a,b,c) \etr\{-(\mathbf{\Omega}_{1} + \mathbf{\Omega}_{2}
    + \mathbf{\Omega}_{3})\}\sum_{\kappa, \tau, \lambda; \phi}^{\infty}\frac{(a + b +c)_{\phi}}{(a)_{\kappa}
    (b)_{\tau}(c)_{\lambda} k! \ t! \ l!}
  \end{equation}
  \par\indent\hfill\mbox{$\displaystyle{\times \ \
    \frac{C_{\phi}^{\kappa,\tau,\lambda}(\mathbf{\Omega}_{1}, \mathbf{\Omega}_{2},
    \mathbf{\Omega}_{3}) C_{\phi}^{\kappa,\tau, \lambda}
  (\mathbf{M}_{1}, \mathbf{M}_{2}, \mathbf{M})} {C_{\phi}(\mathbf{I}_{m})}}$,}
  \par\noindent
and then the nonsymetrised density function of  $\mathbb{U} = (\mathbf{U}_{1}\vdots
\mathbf{U}_{2})' \in \Re^{2m \times m}$ is
\begin{equation}\label{dnbgbt1}\hspace{-0.7cm}
    \mathcal{BGB}I_{2m \times m}(\mathbb{U};a,b,c) \etr\{-(\mathbf{\Omega}_{1} + \mathbf{\Omega}_{2}
    + \mathbf{\Omega}_{3})\}\sum_{\kappa, \tau, \lambda; \phi}^{\infty}\frac{(a + b +c)_{\phi}
    \ \theta_{\phi}^{\kappa,\tau,\lambda}}{(a)_{\kappa} (b)_{\tau}(c)_{\lambda} k! \ t! \ l!}
  \end{equation}
  \par\indent\hfill\mbox{$\displaystyle{\times \ \
   C_{\phi}^{\kappa,\tau, \lambda}(\mathbf{\Omega}_{1} \mathbf{M}_{1}, \mathbf{\Omega}_{2} \mathbf{M}_{2},
   \mathbf{\Omega}_{3} \mathbf{M})}$,}
  \par\noindent
which is denoted as $\mathbb{U} \sim \mathcal{BGB}I_{2m \times m}(a,b,c,
\mathbf{\Omega}_{1}, \mathbf{\Omega}_{2}, \mathbf{\Omega}_{3})$; where $\mathbf{0} <
\mathbf{U}_{1} <\mathbf{I}_{m}$ and $\mathbf{0} < \mathbf{U}_{2} < \mathbf{I}_{m}$
and
\begin{eqnarray*}
  \mathbf{M}_{1} &=& (\mathbf{I}_{m}-\mathbf{U}_{2})(\mathbf{I}_{m}-\mathbf{U}_{1} \mathbf{U}_{2})^{-1}\mathbf{U}_{1}, \\
  \mathbf{M}_{2} &=& (\mathbf{I}_{m}-\mathbf{U}_{1})(\mathbf{I}_{m}-\mathbf{U}_{1} \mathbf{U}_{2})^{-1}\mathbf{U}_{2}, \\
  \mathbf{M} &=& (\mathbf{I}_{m}-\mathbf{U}_{1})(\mathbf{I}_{m}-\mathbf{U}_{1}
  \mathbf{U}_{2})^{-1}(\mathbf{I}_{m}-\mathbf{U}_{2}).
\end{eqnarray*}
\end{thm}
\textit{Proof.} The joint density function of $\mathbf{A}$, $\mathbf{B}$ and
$\mathbf{C}$ is
$$
  \frac{\etr\{-(\mathbf{\Omega}_{1} + \mathbf{\Omega}_{2} + \mathbf{\Omega}_{3})\}}{\Gamma_{m}[a]
  \Gamma_{m}[b]\Gamma_{m}[c]} |\mathbf{A}|^{a-(m+1)/2} |\mathbf{B}|^{b-(m+1)/2} |\mathbf{C}|^{c-(m+1)/2}
  \etr\{-(\mathbf{A} + \mathbf{B} + \mathbf{C})\}
$$
\par\indent\hfill\mbox{$\displaystyle{\times \
    {}_{0}F_{1}(a;\mathbf{\Omega}_{1}\mathbf{A})
  {}_{0}F_{1}(b;\mathbf{\Omega}_{2}\mathbf{B}) {}_{0}F_{1}(c;\mathbf{\Omega}_{3}\mathbf{C})}$,}
  \par\noindent
Now, consider the transformations (\ref{bgb1}) and $\mathbf{C} = \mathbf{C}$, then
$$
  (d\mathbf{A})(d\mathbf{B})(d\mathbf{C}) = |\mathbf{C}|^{m+1} |\mathbf{I}_{m} - \mathbf{U}_{1}|
  ^{-(m+1)} |\mathbf{I}_{m} - \mathbf{U}_{2}|^{-(m+1)}
  (d\mathbf{U}_{1})(d\mathbf{U}_{2})(d\mathbf{C}).
$$
The joint density function of $\mathbf{U}_{1}$, $\mathbf{U}_{2}$ and $\mathbf{C}$
is
\begin{equation}\label{jointUC}
  \frac{\etr\{-(\mathbf{\Omega}_{1} + \mathbf{\Omega}_{2} + \mathbf{\Omega}_{3})\}}
  {\Gamma_{m}[a]\Gamma_{m}[b]\Gamma_{m}[c]} \frac{|\mathbf{U}_{1}|^{a-(m+1)/2} |\mathbf{U}_{2}|^{b-(m+1)/2}}
  {|\mathbf{I}_{m}-\mathbf{U}_{1}|^{a+(m+1)/2}
  |\mathbf{I}_{m}-\mathbf{U}_{2}|^{b+(m+1)/2}} \hspace{1cm}
\end{equation}
\par\indent\hfill\mbox{$\displaystyle{\times \ |\mathbf{C}|^{a+b+c-(m+1)/2}
    \exp\{-(\mathbf{I}_{m}-\mathbf{U}_{1})^{-1} (\mathbf{I}_{m}-\mathbf{U}_{1}\mathbf{U}_{2})
    (\mathbf{I}_{m}-\mathbf{U}_{2})^{-1}\mathbf{C}\},\hspace{1.2cm}}$}
\par\noindent
\par\indent\hfill\mbox{$\displaystyle{\times \
   {}_{0}F_{1}\left(a; \mathbf{\Omega}_{1} \mathbf{C}^{1/2} (\mathbf{I}_{m}-\mathbf{U}_{1})^{-1}
   \mathbf{U}_{1} \mathbf{C}^{1/2}\right)
   {}_{0}F_{1}\left(b; \mathbf{\Omega}_{2} \mathbf{C}^{1/2} (\mathbf{I}_{m}-\mathbf{U}_{2})^{-1}
   \mathbf{U}_{2} \mathbf{C}^{1/2}\right)}$,}
\par\noindent
\par\indent\hfill\mbox{$\displaystyle{\times \
   {}_{0}F_{1}\left(c; \mathbf{\Omega}_{3} \mathbf{C}\right)}$.}
\par\noindent
Now, note that
\begin{description}
\item[i)]
$$
  (\mathbf{I}_{m}-\mathbf{U}_{1})^{-1} - \mathbf{I}_{m}
   =\left\{
        \begin{array}{l}
            (\mathbf{I}_{m}-\mathbf{U}_{1})^{-1}\mathbf{U}_{1}, \\
            \mathbf{U}_{1}(\mathbf{I}_{m}-\mathbf{U}_{1})^{-1}.
        \end{array}
    \right.
$$
With similarly expressions for $(\mathbf{I}_{m}-\mathbf{U}_{2})^{-1} -
\mathbf{I}_{m}$.
\item[ii)] From the argument of $\etr(\cdot)$  in (\ref{jointUC}), denotes $\mathbf{M}^{-1}
= (\mathbf{I}_{m}-\mathbf{U}_{1})^{-1} \mathbf{U}_{1} +
(\mathbf{I}_{m}-\mathbf{U}_{2})^{-1}\mathbf{U}_{2} + \mathbf{I}_{m}$ and note that
$$
\mathbf{M}^{-1} =%
\left\{
  \begin{array}{l}
    (\mathbf{I}_{m}-\mathbf{U}_{1})^{-1}(\mathbf{I}_{m}-\mathbf{U}_{1} \mathbf{U}_{2})
    (\mathbf{I}_{m}-\mathbf{U}_{2})^{-1}, \\
    (\mathbf{I}_{m}-\mathbf{U}_{2})^{-1}(\mathbf{I}_{m}-\mathbf{U}_{1} \mathbf{U}_{2})
    (\mathbf{I}_{m}-\mathbf{U}_{1})^{-1}.
  \end{array}
\right.
$$
\item[iii)] Now assuming that
$(\mathbf{I}_{m}-\mathbf{U}_{1})^{-1}\mathbf{U}_{1}\mathbf{M}$ is an argument of a
symmetric function ($f(\mathbf{AB}) = f(\mathbf{BA})$) by \textbf{i)} and
\textbf{ii)} we have
$$
(\mathbf{I}_{m}-\mathbf{U}_{1})^{-1}\mathbf{U}_{1}\mathbf{M} = %
\left\{
  \begin{array}{l}
    (\mathbf{I}_{m}-\mathbf{U}_{2})(\mathbf{I}_{m}-\mathbf{U}_{1} \mathbf{U}_{2})^{-1}
    \mathbf{U}_{1}\\
    (\mathbf{I}_{m}-\mathbf{U}_{1} \mathbf{U}_{2})^{-1}(\mathbf{I}_{m}-\mathbf{U}_{2})
    \mathbf{U}_{1}, \\
    \mathbf{U}_{1}(\mathbf{I}_{m}-\mathbf{U}_{2})(\mathbf{I}_{m}-\mathbf{U}_{1}
    \mathbf{U}_{2})^{-1},\\
    \mathbf{U}_{1}(\mathbf{I}_{m}-\mathbf{U}_{1} \mathbf{U}_{2})^{-1}(\mathbf{I}_{m} -
    \mathbf{U}_{2}).
  \end{array}
\right.
$$
With similar expressions for
$(\mathbf{I}_{m}-\mathbf{U}_{2})^{-1}\mathbf{U}_{2}\mathbf{M}$.
\end{description}
The marginal joint density function of $\mathbf{U}_{1}$ and $\mathbf{U}_{2}$ is
$$
  \frac{ \etr\{-(\mathbf{\Omega}_{1} + \mathbf{\Omega}_{2}
    + \mathbf{\Omega}_{3})\} |\mathbf{U}_{1}|^{a-(m+1)/2} |\mathbf{U}_{2}|^{b-(m+1)/2}}
    {\Gamma_{m}[a]\Gamma_{m}[b]\Gamma_{m}[c]|\mathbf{I}_{m}-\mathbf{U}_{1}|^{a+(m+1)/2}
    |\mathbf{I}_{m}-\mathbf{U}_{2}|^{b+(m+1)/2}}\hspace{4cm}
$$
\par\indent\hfill\mbox{$\displaystyle{\times \ \
   \int_{\mathbf{C}>\mathbf{0}} |\mathbf{C}|^{a+b+c-(m+1)/2} \exp\{-\mathbf{M^{-1}C}\}
   {}_{0}F_{1}\left(a; \mathbf{\Omega}_{1} \mathbf{C}^{1/2} (\mathbf{I}_{m}-\mathbf{U}_{1})^{-1}
   \mathbf{U}_{1} \mathbf{C}^{1/2}\right) }$,}
\par\noindent
\par\indent\hfill\mbox{$\displaystyle{\times \ \
   {}_{0}F_{1}\left(b; \mathbf{\Omega}_{2} \mathbf{C}^{1/2} (\mathbf{I}_{m}-\mathbf{U}_{2})^{-1}
   \mathbf{U}_{2} \mathbf{C}^{1/2}\right) {}_{0}F_{1}\left(c; \mathbf{\Omega}_{3} \mathbf{C}\right)(d\mathbf{C})}$,}
\par\noindent
Denoting the density joint function of $\mathbf{U}_{1}$ and $\mathbf{U}_{2}$ as
$f_{_{\mathbf{U}_{1},\mathbf{U}_{2}}}(\mathbf{U}_{1},\mathbf{U}_{2})$, considering
the corresponding symmetrised function
$$
  f_{s}(\mathbf{U}_{1},\mathbf{U}_{2}) = \int_{\mathcal{O}(m)}f_{_{\mathbf{U}_{1},\mathbf{U}_{2}}}
  (\mathbf{HU}_{1}\mathbf{H}', \mathbf{HU}_{2}\mathbf{H}')(d\mathbf{H})
$$
and the transformation $\mathbf{C} = \mathbf{HCH}'$ with $(d\mathbf{C}) =
(d(\mathbf{HCH}'))$. Then expanding the hypergeometric functions ${}_{0}F_{1}(\cdot)$
in terms of zonal polynomials we obtain that $f_{s}(\mathbf{U}_{1},\mathbf{U}_{2})$
is
$$
  \frac{ \etr\{-(\mathbf{\Omega}_{1} + \mathbf{\Omega}_{2}
    + \mathbf{\Omega}_{3})\} |\mathbf{U}_{1}|^{a-(m+1)/2} |\mathbf{U}_{2}|^{b-(m+1)/2}}
    {\Gamma_{m}[a]\Gamma_{m}[b]\Gamma_{m}[c]|\mathbf{I}_{m}-\mathbf{U}_{1}|^{a+(m+1)/2}
    |\mathbf{I}_{m}-\mathbf{U}_{2}|^{b+(m+1)/2}}\hspace{4cm}
$$
\par\indent\hfill\mbox{$\displaystyle{\times \ \
   \sum_{k=0}^{\infty} \sum_{t=0}^{\infty} \sum_{l=0}^{\infty}  \sum_{\kappa}  \sum_{\tau}
   \sum_{\lambda} \frac{1}{(a)_{\kappa} (b)_{\tau} (c)_{\lambda} k! t! l!}
   \int_{\mathbf{C}>\mathbf{0}} |\mathbf{C}|^{a+b+c-(m+1)/2} \exp\{-\mathbf{M^{-1}C}\}
    }$,}
\par\noindent
\par\indent\hfill\mbox{$\displaystyle{\times \ \
   \left [ \int_{\mathcal{O}(m)} C_{\kappa}\left(\mathbf{\Omega}_{1} \mathbf{HC}^{1/2}
   (\mathbf{I}_{m}-\mathbf{U}_{1})^{-1}\mathbf{U}_{1} \mathbf{C}^{1/2}\mathbf{H}'\right) \right .}$,}
\par\noindent
\par\indent\hfill\mbox{$\displaystyle{\times \ \ \left.
   C_{\tau}\left(\mathbf{\Omega}_{2} \mathbf{HC}^{1/2} (\mathbf{I}_{m}-\mathbf{U}_{2})^{-1}
   \mathbf{U}_{2} \mathbf{C}^{1/2}\mathbf{H}'\right) C_{\lambda}\left(\mathbf{\Omega}_{3}
   \mathbf{HC}\mathbf{H}'\right) (d\mathbf{H})\right ](d\mathbf{C})}$,}
\par\noindent
By integrating with respect to $\mathbf{H}$, using \citet[equation (2.2)]{chd:86} and
the notation for the operator sum as in \citet{da:80}, we have
$$
  \frac{ \etr\{-(\mathbf{\Omega}_{1} + \mathbf{\Omega}_{2}
    + \mathbf{\Omega}_{3})\} |\mathbf{U}_{1}|^{a-(m+1)/2} |\mathbf{U}_{2}|^{b-(m+1)/2}}
    {\Gamma_{m}[a]\Gamma_{m}[b]\Gamma_{m}[c]|\mathbf{I}_{m}-\mathbf{U}_{1}|^{a+(m+1)/2}
    |\mathbf{I}_{m}-\mathbf{U}_{2}|^{b+(m+1)/2}}\hspace{3.5cm}
$$
\par\indent\hfill\mbox{$\displaystyle{\times \ \
   \sum_{\kappa, \tau, \lambda; \phi}^{\infty} \frac{C_{\phi}^{\kappa,\tau,\lambda}
   (\mathbf{\Omega}_{1},\mathbf{\Omega}_{2}, \mathbf{\Omega}_{3})}{(a)_{\kappa} (b)_{\tau} (c)_{\lambda}
   k! t! l!}\int_{\mathbf{C}>\mathbf{0}} |\mathbf{C}|^{a+b+c-(m+1)/2} \exp\{-\mathbf{M^{-1}C}\}}$,}
\par\noindent
\par\indent\hfill\mbox{$\displaystyle{\times \ \
    \frac{C_{\phi}^{\kappa,\tau,\lambda}\left((\mathbf{I}_{m}-\mathbf{U}_{1})^{-1}\mathbf{U}_{1}\mathbf{M}^{-1},
    (\mathbf{I}_{m}-\mathbf{U}_{2})^{-1} \mathbf{U}_{2}\mathbf{M}^{-1}, \mathbf{M}^{-1}\right)}
    {C_{\phi}(\mathbf{I}_{m})}(d\mathbf{C})}$}.
\par\noindent
Finality, by integrating with respect to $\mathbf{C}$, see \citet[equation
(3.21)]{ch:80} and \textbf{iii)}, we obtain the joint symmetrised density function of
$\mathbf{U}_{1}$ and $\mathbf{U}_{2}$.

The joint nonsymmetrised density function of $\mathbf{U}_{1}$ and $\mathbf{U}_{2}$ is
obtained by applying the idea of \citet{g:73} in an inverse way. With this proposal
observe that $\mathcal{BGB}I_{2m \times m}(\mathbf{H}\mathbb{U}\mathbf{H}';a,b,c) =
\mathcal{BGB}I_{2m \times m}(\mathbb{U};a,b,c)$ and by \citet{dg:06}
$$
  \int_{\mathcal{O}(m)}C_{\phi}^{\kappa,\tau, \lambda}(\mathbf{\Omega}_{1} \mathbf{H}
  \mathbf{M}_{1}\mathbf{H}', \mathbf{\Omega}_{1} \mathbf{H}\mathbf{M}_{2}\mathbf{H}',
   \mathbf{\Omega}_{1} \mathbf{H} \mathbf{M} \mathbf{H}') (d\mathbf{H})\hspace{5cm}
$$
\par\indent\hfill\mbox{$\displaystyle{
    = \frac{C_{\phi}^{\kappa,\tau,\lambda}(\mathbf{\Omega}_{1}, \mathbf{\Omega}_{2},
    \mathbf{\Omega}_{3})
   C_{\phi}^{\kappa,\tau, \lambda}(\mathbf{M}_{1}, \mathbf{M}_{2}, \mathbf{M})}
   {\theta_{\phi}^{\kappa,\tau, \lambda}C_{\phi}(\mathbf{I}_{m})}}$},
\par\noindent
from where the desired result is obtained. \qed

In addition, note that in Theorem \ref{dngb1},
$$
  \mathbf{U}_{1} \sim \mathcal{B}I_{m}(a, c, {\mathbf\Omega}_{1}, \mathbf{\Omega}_{3})
  \quad \mbox{and} \quad \mathbf{U}_{2}\sim \mathcal{B}I_{m}(b, c, \mathbf{\Omega}_{2},
  \mathbf{\Omega}_{3}).
$$

Next, assuming that one and/or two of the matrices $\mathbf{A}$, $\mathbf{B}$ or
$\mathbf{C}$ have a central matrix variate gamma distribution in Theorem \ref{dngb1},
let us study all the possible nonsymetrised noncentral densities.

\begin{cor} \label{coro1}Let us assume the hypothesis of Theorem \ref{dngb1}. Then the joint nonsymetrised density
function of $\mathbf{U}_{1}$ and $\mathbf{U}_{2}$ is:
\begin{enumerate}
  \item If $\mathbf{\Omega}_{1}= \mathbf{\Omega}_{2} = \mathbf{0}$
  $$
    \frac{\mathcal{BGB}I_{2m \times
    m}(\mathbb{U};a,b,c)}{\etr\{\mathbf{\Omega}_{3}\}}
    {}_{1}F_{1}(a+b+c; c; \mathbf{\Omega}_{3} \mathbf{M})
  $$
  and $\mathbf{U}_{1} \sim \mathcal{B}I(A)_{m}(a, c, \mathbf{\Omega}_{3})$
  and $\mathbf{U}_{2} \sim \mathcal{B}I(A)_{m}(b, c, \mathbf{\Omega}_{3})$.
  \item If $\mathbf{\Omega}_{3}= \mathbf{0}$
  $$
    \frac{\mathcal{BGB}I_{2m \times m}(\mathbb{U};a,b,c)}{\etr\{\mathbf{\Omega}_{1} +
    \mathbf{\Omega}_{2}\}}\sum_{\kappa, \tau; \phi}^{\infty}\frac{(a + b +c)_{\phi}
    \ \theta_{\phi}^{\kappa,\tau}}{(a)_{\kappa} (b)_{\tau} k! \ t!}
    C_{\phi}^{\kappa,\tau}(\mathbf{\Omega}_{1} \mathbf{M}_{1}, \mathbf{\Omega}_{2} \mathbf{M}_{2})
  $$
  also we have that $\mathbf{U}_{1} \sim \mathcal{B}I(B)_{m}(a, c, \mathbf{\Omega}_{1})$
  and $\mathbf{U}_{2} \sim \mathcal{B}I(B)_{m}(b, c, \mathbf{\Omega}_{2})$.
  \item If $\mathbf{\Omega}_{2}= \mathbf{\Omega}_{3} = \mathbf{0}$
  $$
    \frac{\mathcal{BGB}I_{2m \times
    m}(\mathbb{U};a,b,c)}{\etr\{\mathbf{\Omega}_{1}\}}
    {}_{1}F_{1}(a+b+c; c; \mathbf{\Omega}_{1} \mathbf{M}_{1})
  $$
  and with $\mathbf{U}_{1} \sim \mathcal{B}I(B)_{m}(a, c, \mathbf{\Omega}_{1})$ and
  $\mathbf{U}_{2} \sim \mathcal{B}I_{m}(b, c)$.
  \item If $\mathbf{\Omega}_{1}= \mathbf{\Omega}_{3} = \mathbf{0}$
  $$
    \frac{\mathcal{BGB}I_{2m \times
    m}(\mathbb{U};a,b,c)}{\etr\{\mathbf{\Omega}_{2}\}}
    {}_{1}F_{1}(a+b+c; c; \mathbf{\Omega}_{2} \mathbf{M}_{2})
  $$
  and $\mathbf{U}_{1} \sim \mathcal{B}I_{m}(a, c)$ and $\mathbf{U}_{2}
  \sim \mathcal{B}I(B)_{m}(b, c, \mathbf{\Omega}_{2})$.
   \item If $\mathbf{\Omega}_{2}= \mathbf{0}$
  $$
    \frac{\mathcal{BGB}I_{2m \times m}(\mathbb{U};a,b,c)}{\etr\{\mathbf{\Omega}_{1} +
    \mathbf{\Omega}_{3}\}}\sum_{\kappa, \lambda; \phi}^{\infty}\frac{(a + b +c)_{\phi}
    \ \theta_{\phi}^{\kappa,\lambda}}{(a)_{\kappa} (c)_{\lambda} k! \ l!}
    C_{\phi}^{\kappa,\lambda}(\mathbf{\Omega}_{1} \mathbf{M}_{1}, \mathbf{\Omega}_{3}
    \mathbf{M}).
  $$
  In addition, note that $\mathbf{U}_{1} \sim \mathcal{B}I_{m}(a, c, \mathbf{\Omega}_{1},
  \mathbf{\Omega}_{3})$ and $\mathbf{U}_{2} \sim \mathcal{B}I(A)_{m}(b, c,\mathbf{\Omega}_{3})$.
   \item If $\mathbf{\Omega}_{1}= \mathbf{0}$
  $$
    \frac{\mathcal{BGB}I_{2m \times m}(\mathbb{U};a,b,c)}{\etr\{\mathbf{\Omega}_{2} +
    \mathbf{\Omega}_{3}\}}\sum_{\tau, \lambda; \phi}^{\infty}\frac{(a + b +c)_{\phi}
    \ \theta_{\phi}^{\tau,\lambda}}{(b)_{\tau} (c)_{\lambda} t! \ l!}
    C_{\phi}^{\tau,\lambda}(\mathbf{\Omega}_{2} \mathbf{M}_{2}, \mathbf{\Omega}_{3} \mathbf{M})
  $$
  and where $\mathbf{U}_{1} \sim \mathcal{B}I(A)_{m}(a, c, \mathbf{\Omega}_{3})$ and $\mathbf{U}_{2}
  \sim \mathcal{B}I_{m}(b, c, \mathbf{\Omega}_{2}, \mathbf{\Omega}_{3})$.
\end{enumerate}
\end{cor}
\textbf{Proof.} The joint density functions of $\mathbf{U}_{1}$ and $\mathbf{U}_{2}$
in all items are a consequence of the basic properties of invariant polynomials, see
\citet[equations (2.1) and (2.3)]{da:79}, see also \citet[equations (3.3) and
(3.6)]{ch:80}. The second claim in each of the items is a consequence of construction
(\ref{bgb12}). \qed

\section{Doubly noncentral bimatrix variate generalised beta type II distribution}\label{sec4}

In this section we derive the doubly noncentral and noncentral bimatrix variate
generalised beta type II distributions.

\begin{thm}\label{dngb2}
Let $\mathbf{A}$, $\mathbf{B}$ and $\mathbf{C}$ be independent random matrices, such
that $\mathbf{A} \sim \mathcal{G}_{m}(a, \mathbf{I}_{m}, \mathbf{\Omega}_{1})$,
$\mathbf{B} \sim \mathcal{G}_{m}(b, \mathbf{I}_{m}, \mathbf{\Omega}_{2})$ and
$\mathbf{C} \sim \mathcal{G}_{m}(c, \mathbf{I}_{m}, \mathbf{\Omega}_{3})$ with
$\R(a)> (m-1)/2$, $\R(b) > (m-1)/2$ and $\R(c) > (m-1)/2$ and let us define
\begin{equation}\label{bgb22}\hspace{-.7cm}
    \mathbf{F}_{1} = \mathbf{C}^{-1/2}\mathbf{A}\mathbf{C}^{-1/2}
  \quad \mbox{and} \quad \mathbf{F}_{2} = \mathbf{C}^{-1/2}\mathbf{B}
  \mathbf{C}^{-1/2}
\end{equation}
Then the symmetrised density function of  $\mathbb{F} = (\mathbf{F}_{1}\vdots
\mathbf{F}_{2})' \in \Re^{2m \times m}$ is
\begin{equation}\label{sdnbgbt2}\hspace{-0.7cm}
    \mathcal{BGB}II_{2m \times m}(\mathbb{F};a,b,c) \etr\{-(\mathbf{\Omega}_{1} + \mathbf{\Omega}_{2}
    + \mathbf{\Omega}_{3})\}\sum_{\kappa, \tau, \lambda; \phi}^{\infty}\frac{(a + b +c)_{\phi}}{(a)_{\kappa}
    (b)_{\tau}(c)_{\lambda} k! \ t! \ l!}
  \end{equation}
  \par\indent\hfill\mbox{$\displaystyle{\times \ \
    \frac{C_{\phi}^{\kappa,\tau,\lambda}(\mathbf{\Omega}_{1}, \mathbf{\Omega}_{2},
    \mathbf{\Omega}_{3}) C_{\phi}^{\kappa,\tau, \lambda}
  (\mathbf{N}_{1}, \mathbf{N}_{2}, \mathbf{N})} {C_{\phi}(\mathbf{I}_{m})}}$,}
  \par\noindent
and then the nonsymetrised density function of  $\mathbb{F} = (\mathbf{F}_{1}\vdots
\mathbf{F}_{2})' \in \Re^{2m \times m}$ is
\begin{equation}\label{dnbgbt2}\hspace{-0.7cm}
    \mathcal{BGB}II_{2m \times m}(\mathbb{F};a,b,c) \etr\{-(\mathbf{\Omega}_{1} + \mathbf{\Omega}_{2}
    + \mathbf{\Omega}_{3})\}\sum_{\kappa, \tau, \lambda; \phi}^{\infty}\frac{(a + b +c)_{\phi}
    \ \theta_{\phi}^{\kappa,\tau,\lambda}}{(a)_{\kappa} (b)_{\tau}(c)_{\lambda} k! \ t! \ l!}
  \end{equation}
  \par\indent\hfill\mbox{$\displaystyle{\times \ \
   C_{\phi}^{\kappa,\tau, \lambda}(\mathbf{\Omega}_{1} \mathbf{N}_{1}, \mathbf{\Omega}_{2} \mathbf{N}_{2},
   \mathbf{\Omega}_{3} \mathbf{N})}$,}
  \par\noindent
which is denoted as $\mathbb{F} \sim \mathcal{BGB}II_{2m \times m}(a,b,c,
\mathbf{\Omega}_{1}, \mathbf{\Omega}_{2}, \mathbf{\Omega}_{3})$; where $
\mathbf{F}_{1} > \mathbf{0}$ and $\mathbf{F}_{2} > \mathbf{0}$ and
\begin{eqnarray*}
  \mathbf{N}_{1} &=& (\mathbf{I}_{m}+\mathbf{F}_{1})^{-1}(\mathbf{I}_{m}+\mathbf{F}_{2})
  (\mathbf{I}_{m}+\mathbf{F}_{1} + \mathbf{F}_{2})^{-1}\mathbf{F}_{1}, \\
  \mathbf{N}_{2} &=& (\mathbf{I}_{m}+\mathbf{F}_{1} + \mathbf{F}_{2})^{-1}
  (\mathbf{I}_{m}+\mathbf{F}_{1})(\mathbf{I}_{m}+\mathbf{F}_{2})^{-1}\mathbf{F}_{2}, \\
  \mathbf{N} &=& (\mathbf{I}_{m}+\mathbf{F}_{1}+ \mathbf{F}_{2})^{-1}.
\end{eqnarray*}
\end{thm}
\textit{Proof.} The proof is similar to that given for Theorem \ref{dngb2}. Or
alternatively, by noting that if $\mathbf{U} \sim \mathcal{B}I_{m}(a,b)$, then
$(\mathbf{I}_{m}-\mathbf{U})^{-1} -\mathbf{I}_{m} \sim \mathcal{B}II_{m}(a,b)$, see
\citet{sk:79} and \citet{dggj:07}. Then, by construction (\ref{bgb12}) we have that
if $\mathbb{U} \sim \mathcal{BGB}I_{2m \times m}(a,b,c, \mathbf{\Omega}_{1},
\mathbf{\Omega}_{2}, \mathbf{\Omega}_{3})$ then
$\mathbb{F}=(\mathbf{F}_{1}|\mathbf{F}_{2})' \sim \mathcal{BGB}II_{2m \times
m}(a,b,c, \mathbf{\Omega}_{1}, \mathbf{\Omega}_{2}, \mathbf{\Omega}_{3})$, where
$$
  \mathbb{F} =
  \left (
    \begin{array}{c}
      \mathbf{F}_{1} \\
      \mathbf{F}_{2}
    \end{array}
  \right ) =
  \left (
    \begin{array}{c}
      (\mathbf{I}_{m}-\mathbf{U}_{1})^{-1}-\mathbf{I}_{m} \\
      (\mathbf{I}_{m}-\mathbf{U}_{2})^{-1}-\mathbf{I}_{m}
    \end{array}
  \right ).
$$
Then the inverse transformation is given by
\begin{equation}\label{itrans}
  \mathbb{U} =
  \left (
    \begin{array}{c}
      \mathbf{U}_{1} \\
      \mathbf{U}_{2}
    \end{array}
  \right ) =
  \left (
    \begin{array}{c}
       \mathbf{I}_{m} - (\mathbf{I}_{m}+\mathbf{F}_{1})^{-1} \\
       \mathbf{I}_{m} - (\mathbf{I}_{m}+\mathbf{F}_{2})^{-1}
    \end{array}
  \right ),
\end{equation}
and the Jacobian is given by
$$
  (d\mathbf{U}_{1})(d\mathbf{U}_{2}) =  |\mathbf{I}_{m} + \mathbf{F}_{1}|^{-(m+1)}
  |\mathbf{I}_{m} + \mathbf{F}_{2}|^{-(m+1)} (d\mathbf{F}_{1})(d\mathbf{F}_{2}).
$$
Then, (\ref{sdnbgbt2}) follows, observing that under transformation (\ref{itrans})
\begin{description}
  \item[i)]
  $$
    \frac{\mathcal{BGB}I_{2m \times m}((\mathbf{I}_{m} + (\mathbf{I}_{m} + \mathbf{F}_{1})^{-1}|
    (\mathbf{I}_{m} - (\mathbf{I}_{m} + \mathbf{F}_{2})^{-1})';a,b,c)}{
    |\mathbf{I}_{m} - \mathbf{F}_{1}|^{(m+1)} |\mathbf{I}_{m} +
    \mathbf{F}_{2}|^{(m+1)}}
  $$
  \par\indent\hfill\mbox{$\displaystyle{ =  \mathcal{BGB}II_{2m \times m}(\mathbb{F};a,b,c)}$,}
  \par\noindent
  \item[ii)] and
  $$
  \begin{array}{lllll}
    \mathbf{M}_{1} &=& (\mathbf{I}_{m}+\mathbf{F}_{1})^{-1}(\mathbf{I}_{m}+\mathbf{F}_{2})
  (\mathbf{I}_{m}+\mathbf{F}_{1} + \mathbf{F}_{2})^{-1}\mathbf{F}_{1} &=& \mathbf{N}_{1}, \\
  \mathbf{M}_{2} &=& (\mathbf{I}_{m}+\mathbf{F}_{1} + \mathbf{F}_{2})^{-1}
  (\mathbf{I}_{m}+\mathbf{F}_{1})(\mathbf{I}_{m}+\mathbf{F}_{2})^{-1}\mathbf{F}_{2} &=& \mathbf{N}_{2}, \\
  \mathbf{M} &=& (\mathbf{I}_{m}+\mathbf{F}_{1}+ \mathbf{F}_{2})^{-1} &=& \mathbf{N}.
  \end{array}
  $$
To obtain (\ref{dnbgbt2}), observe that
$$
  \mathcal{BGB}II_{2m \times m}(\mathbf{H}\mathbb{F}\mathbf{H}';a,b,c) = \mathcal{BGB}II_{2m \times
   m}(\mathbb{F};a,b,c)
$$
and by \citet{dg:06}
$$
  \int_{\mathcal{O}(m)}C_{\phi}^{\kappa,\tau, \lambda}(\mathbf{\Omega}_{1} \mathbf{H}
  \mathbf{N}_{1}\mathbf{H}', \mathbf{\Omega}_{1} \mathbf{H}\mathbf{N}_{2}\mathbf{H}',
   \mathbf{\Omega}_{1} \mathbf{H} \mathbf{N} \mathbf{H}') (d\mathbf{H})\hspace{5cm}
$$
\par\indent\hfill\mbox{$\displaystyle{
    = \frac{C_{\phi}^{\kappa,\tau,\lambda}(\mathbf{\Omega}_{1}, \mathbf{\Omega}_{2}, \mathbf{\Omega}_{3})
   C_{\phi}^{\kappa,\tau, \lambda}(\mathbf{N}_{1}, \mathbf{N}_{2}, \mathbf{N})}
   {\theta_{\phi}^{\kappa,\tau, \lambda}C_{\phi}(\mathbf{I}_{m})}}$},
\par\noindent
from where the nonsymmetrised joint density function of $\mathbf{F}_{1}$ and
$\mathbf{F}_{2}$ is obtained. \qed

Finally, let us find the different possibilities of the nonsymetrised noncentral
density functions of bimatrix variate generalised beta type II distributions.

\begin{cor} \label{coro2} Under the hypothesis of Theorem \ref{dngb2} the joint nonsymetrised density
function of $\mathbf{F}_{1}$ and $\mathbf{F}_{2}$ is:
\begin{enumerate}
  \item If $\mathbf{\Omega}_{1}= \mathbf{\Omega}_{2} = \mathbf{0}$
  $$
    \frac{\mathcal{BGB}II_{2m \times m}(\mathbb{F};a,b,c)}{\etr\{\mathbf{\Omega}_{3}\}}
    {}_{1}F_{1}(a+b+c; c; \mathbf{\Omega}_{3} \mathbf{N})
  $$
  and $\mathbf{F}_{1} \sim \mathcal{B}II(A)_{m}(a, c, \mathbf{\Omega}_{3})$
  and $\mathbf{F}_{2} \sim \mathcal{B}II(A)_{m}(b, c, \mathbf{\Omega}_{3})$.
  \item If $\mathbf{\Omega}_{3}= \mathbf{0}$
  $$
    \frac{\mathcal{BGB}II_{2m \times m}(\mathbb{F};a,b,c)}{\etr\{\mathbf{\Omega}_{1} +
    \mathbf{\Omega}_{2}\}}\sum_{\kappa, \tau; \phi}^{\infty}\frac{(a + b +c)_{\phi}
    \ \theta_{\phi}^{\kappa,\tau}}{(a)_{\kappa} (b)_{\tau} k! \ t!}
    C_{\phi}^{\kappa,\tau}(\mathbf{\Omega}_{1} \mathbf{N}_{1}, \mathbf{\Omega}_{2}
    \mathbf{N}_{2}).
  $$
  In addition,  we have that $\mathbf{F}_{1} \sim \mathcal{B}II(B)_{m}(a, c, \mathbf{\Omega}_{1})$
  and $\mathbf{F}_{2} \sim \mathcal{B}II(B)_{m}(b, c, \mathbf{\Omega}_{2})$.
  \item If $\mathbf{\Omega}_{2}= \mathbf{\Omega}_{3} = \mathbf{0}$
  $$
    \frac{\mathcal{BGB}II_{2m \times
    m}(\mathbb{F};a,b,c)}{\etr\{\mathbf{\Omega}_{1}\}}
    {}_{1}F_{1}(a+b+c; c; \mathbf{\Omega}_{1} \mathbf{N}_{1})
  $$
  and with $\mathbf{F}_{1} \sim \mathcal{B}II(B)_{m}(a, c, \mathbf{\Omega}_{1})$ and
  $\mathbf{F}_{2} \sim \mathcal{B}II_{m}(b, c)$.
  \item If $\mathbf{\Omega}_{1}= \mathbf{\Omega}_{3} = \mathbf{0}$
  $$
    \frac{\mathcal{BGB}II_{2m \times
    m}(\mathbb{F};a,b,c)}{\etr\{\mathbf{\Omega}_{2}\}}
    {}_{1}F_{1}(a+b+c; c; \mathbf{\Omega}_{2} \mathbf{N}_{2})
  $$
  and $\mathbf{F}_{1} \sim \mathcal{B}II_{m}(a, c)$ and $\mathbf{F}_{2}
  \sim \mathcal{B}II(B)_{m}(b, c, \mathbf{\Omega}_{2})$.
   \item If $\mathbf{\Omega}_{2}= \mathbf{0}$
  $$
    \frac{\mathcal{BGB}II_{2m \times m}(\mathbb{F};a,b,c)}{\etr\{\mathbf{\Omega}_{1} +
    \mathbf{\Omega}_{3}\}}\sum_{\kappa, \lambda; \phi}^{\infty}\frac{(a + b +c)_{\phi}
    \ \theta_{\phi}^{\kappa,\lambda}}{(a)_{\kappa} (c)_{\lambda} k! \ l!}
    C_{\phi}^{\kappa,\lambda}(\mathbf{\Omega}_{1} \mathbf{N}_{1}, \mathbf{\Omega}_{3}
    \mathbf{N}).
  $$
  In addition note that $\mathbf{F}_{1} \sim \mathcal{B}II_{m}(a, c, \mathbf{\Omega}_{1},
  \mathbf{\Omega}_{3})$ and $\mathbf{F}_{2} \sim \mathcal{B}II(A)_{m}(b, c,\mathbf{\Omega}_{3})$.
   \item If $\mathbf{\Omega}_{1}= \mathbf{0}$
  $$
    \frac{\mathcal{BGB}II_{2m \times m}(\mathbb{F};a,b,c)}{\etr\{\mathbf{\Omega}_{2} +
    \mathbf{\Omega}_{3}\}}\sum_{\tau, \lambda; \phi}^{\infty}\frac{(a + b +c)_{\phi}
    \ \theta_{\phi}^{\tau,\lambda}}{(b)_{\tau} (c)_{\lambda} t! \ l!}
    C_{\phi}^{\tau,\lambda}(\mathbf{\Omega}_{2} \mathbf{N}_{2}, \mathbf{\Omega}_{3} \mathbf{N})
  $$
  and where $\mathbf{F}_{1} \sim \mathcal{B}II(A)_{m}(a, c, \mathbf{\Omega}_{3})$ and $\mathbf{F}_{2}
  \sim \mathcal{B}II_{m}(b, c, \mathbf{\Omega}_{2}, \mathbf{\Omega}_{3})$.
\end{enumerate}
\end{cor}
\textit{Proof.} The proof is obtained in a similar way to the proof of Corollary
\ref{coro1}. \qed
\end{description}

\section{Conclusions}

The problem to finding the density function of a doubly noncentral beta type II
distribution has been studied by different authors, see \citet{da:79}, \citet{ch:80}
and \citet{dggj:08a}, among others. All these authors, in fact, found the symmetrised
density function of a doubly noncentral beta type II distribution. The nonsymmetrised
density function remained an open problem. The Theorem \ref{teo1} solves this problem
for doubly noncentral beta type I and II distributions by applying in an inverse way
the definition of symmetrised function proposed by \citet{g:73}. This approach was
used previously by \cite{dggj:07} in the noncentral cases.

In a similar way, we have found the symmetrised doubly noncentral generalised beta
type I and II distributions and, by applying in an inverse way the definition of
symmetrised function proposed by \citet{g:73}, we have found the noncymmetrised
doubly noncentral generalised beta type I and II distributions. Finally, as
corollaries, we studied all the different noncentral generalised beta type I and II
distributions.

\section*{Acknowledgements} This research work was partially supported  by
CONACYT-M\'exico research grant No. \ 81512 and IDI-Spain, grants FQM2006-2271 and
MTM2008-05785. This paper was written during J. A. D\'{\i}az- Garc\'{\i}a's stay as a
visiting professor at the Department of Statistics and O. R. of the University of
Granada, Spain.

\end{document}